\documentclass[4pt]{article} 

\usepackage[utf8]{inputenc} 
\usepackage{geometry} 
\usepackage{graphicx} 

\usepackage[colorlinks,citecolor=red]{hyperref}
\usepackage{amsmath}
\usepackage{amsfonts}
\usepackage{dsfont}
\usepackage{mathrsfs}
\usepackage{latexsym}
\usepackage{graphicx}
\usepackage{amssymb}
\usepackage{float}
\usepackage{epsfig}
\usepackage{epstopdf}
\usepackage{color,xcolor}
\usepackage{booktabs} 
\usepackage{array} 
\usepackage{paralist} 
\usepackage{verbatim} 
\usepackage{subfig} 

\newtheorem{theorem}{Theorem}[section]
\newtheorem{lemma}[theorem]{Lemma}
\newtheorem{corollary}[theorem]{Corollary}

\usepackage{fancyhdr} 
\usepackage{sectsty}
\allsectionsfont{\sffamily\mdseries\upshape} 

\usepackage[nottoc,notlof,notlot]{tocbibind} 
\usepackage[titles,subfigure]{tocloft} 

\title{A note on eigenvalues and Hamiltoinan properties of $k$-connected graphs}

\author{Huicai Jia \thanks{College of Science, Henan
University of Engineering, Zhengzhou, Henan 451191, China; School of Mathematics, Renmin University of China, Beijing, 100872, China. Email: jhc607@163.com},
Ruifang Liu \thanks{Corresponding author. School of Mathematics and Statistics, Zhengzhou
University, Zhengzhou, Henan 450001, China. Email: rfliu@zzu.edu.cn},
Hong-Jian Lai \thanks{Department of Mathematics, West Virginia
University, Morgantown, WV 26506, USA. E-mail: hjlai@math.wvu.edu}}
\date{} 

\begin{document}

\maketitle

\begin{abstract}
Let $\lambda_{1}(G)$ and $\mu_{1}(G)$ denote the spectral radius and the Laplacian spectral radius of a graph $G$, respectively.
Li in [Electronic J. Linear Algebra 34 (2018) 389-392] proved sharp upper bounds of $\lambda_{1}(G)$ based on the connectivity to assure a connected graph
to be Hamiltonian and traceable, respectively. In this paper, we present best possible upper bounds of $\lambda_{1}(G)$ for $k$-connected graphs to be
Hamiltonian-connected and homogeneously traceable, respectively. Furthermore, best possible upper bounds of $\mu_{1}(G)$ to predict $k$-connected graphs
to be Hamiltonian-connected, Hamiltonian and traceable are originally proved, respectively.
\end{abstract}



\section{Introduction}
We consider simple, undirected and connected graphs with undefined terms and notation reference to
\cite{BM}. As in \cite{BM}, $\overline{G}$, $\alpha(G)$, $\kappa(G)$, $\delta(G)$ and $d(v)$ denote the
{\bf complement}, the {\bf stability number} (also call the {\bf independence number}),
the connectivity, the minimum degree of a graph $G$ and the degree of $v$ in $G$,
respectively.
The {\bf join} of $G$ and $H$, denoted by $G\nabla H$, is the graph obtained from a disjoint
union of $G$ and $H$ by adding all possible edges between them. Let $K_{a, b}$ denote complete bipartite graphs
on $n$ vertices, where $a+b=n.$

A well-known result of Whitney \cite{Whit32} states that $\kappa(G)\leq \delta(G)$ for any graph $G$.
A graph $G$ is {\bf $k$-connected} if
$\kappa(G)\geq k$.
A cycle passing through all the vertices of a graph is called a Hamiltonian cycle. A path passing through all
the vertices of a graph is called a Hamiltonian path. The graph $G$ is called {\bf Hamiltonian-connected} if every two vertices of $G$ are
connected by a Hamiltonian path. A graph containing a Hamiltonian cycle is called a {\bf Hamiltonian graph}. Surely all Hamilton-connected graphs are Hamiltonian. A graph $G$ is called {\bf homogeneously traceable} if for each $v\in V(G)$, there is a Hamiltonian path in $G$ with initial vertex
$v.$ Clearly every Hamiltonian graph is homogeneously traceable. A graph containing a Hamiltonian
path is said to be {\bf traceable}. Hence these four Hamiltonian properties weaken in turn.

The {\bf adjacency matrix} of $G$ is the
$n\times n$ matrix $A(G)=(a_{ij})$, where $a_{ij}=1$ if $v_i$ and $v_j$ are adjacent and
otherwise $a_{ij}=0$. Let $\lambda_1(G)\geq \lambda_2(G)\geq
\cdots \geq \lambda_n(G)$ be the adjacency spectrum of $G$.
Let $D(G)$ be the diagonal matrix of the vertex degrees of $G$.
The matrix $L(G) = D(G)-A(G)$ is known as the {\bf Laplacian matrix}  of $G$. Let $\mu_1(G)\geq \mu_2(G)\geq
\cdots \geq \mu_n(G)$ be the Laplacian eigenvalues of $G$.
We call $\lambda_1(G)$ and $\mu_1(G)$ the {\bf spectral radius} and the {\bf Laplacian spectral radius} of $G$, respectively.

There have been many studies on spectral
conditions which warrant Hamiltonian properties of a graph, as can be seen in
\cite{FN, PX, LRF, LM, NB, YGD, ZB, ZQ}, among others. However, hardly any of these results involve Laplacian eigenvalues.
Very recently, R. Li \cite{LR} proved sufficient conditions of $\lambda_{1}(G)$ based on the connectivity to
assure a connected graph to be Hamiltonian and traceable, respectively.

\begin{theorem}\label{old} (R. Li \cite{LR}) Let $G$ be a graph of order $n\geq3$ with connectivity $\kappa$ and minimum degree $\delta(G).$
If $\displaystyle \lambda_{1}(G)\leq \delta \sqrt{\frac{\kappa-s+1}{n-\kappa+s-1}},$ then each of the following holds.\\
(i) $G$ is Hamiltonian if and only if $G\ncong K_{\kappa, \kappa-s+1}$ for $s=0$.\\
(ii)  $G$ is traceable if and only if $G\ncong K_{\kappa, \kappa-s+1}$ for $s=-1$, where $n\geq12.$
\end{theorem}

One of our goals is to show best possible upper bounds of $\lambda_{1}(G)$ for $k$-connected graphs to be
Hamiltonian-connected and homogeneously traceable, respectively.
Another goal of this research is to initiate studies to find best possible upper bounds of $\mu_{1}(G)$
to predict $k$-connected graphs to be Hamiltonian-connected, Hamiltonian and traceable, respectively.
The main results are as follows.

\begin{theorem}\label{main1} Let $G$ be a $k$-connected graph of order $n\geq3$ and minimum degree $\delta(G).$
If $\displaystyle \lambda_{1}(G)\leq \delta \sqrt{\frac{k-s+1}{n-k+s-1}}$,
then  $G$ is Hamiltonian-connected if and only if $G\ncong K_{k, k-s+1}$ for $s=1$.
\end{theorem}
As a consequence of Theorem \ref{main1}, an upper bound on $\lambda_{1}(K_{1}\nabla G)$
to assure a $k$-connected graph to be homogeneously traceable is obtained.

\begin{corollary}\label{main2}
Let $G$ be a $k$-connected graph of order $n\geq3$ and minimum degree $\delta(G).$ If   $\displaystyle \lambda_{1}(K_{1}\nabla G)\leq(\delta+1)\sqrt{\frac{k+1}{n-k}},$
then $G$ is homogeneously traceable.
\end{corollary}

\begin{theorem}\label{main3}  Let $G$ be a $k$-connected graph of order $n\geq3$ and minimum degree $\delta(G).$ If $\displaystyle \mu_{1}(G)<\frac{n\delta}{n-k+s-1},$
then each of the following holds.\\
(i) $G$ is Hamiltonian-connected for $s=1$.\\
(ii) $G$ is Hamiltonian for $s=0$.\\
(iii) $G$ is traceable for $s=-1$.
\end{theorem}

Noting that the results of Theorem \ref{main3} are also best possible in the sense that the condition of the theorem can not be weakened.
Let us consider $G\cong K_{k, k-s+1}.$ From Lemma \ref{le6}, we have $\mu_{1}(K_{k, k-s+1})=2k-s+1.$ Clearly the graph $G\cong K_{k, k-s+1}$ satisfies $\mu_{1}(G)=\frac{n\delta}{n-k+s-1}$ for $s=1, 0, -1.$ However, $G\cong K_{k,k}$ is not Hamiltonian-connected for $s=1$,
$G\cong K_{k,k+1}$ is not Hamiltonian for $s=0$. $G\cong K_{k,k+2}$ is not traceable for $s =-1$.


In the next section, we display some tools to be deployed in
our arguments. The proofs of the main results are in the subsequent
section.

\section{Preliminaries}

We in this section will present some lemmas that will be useful in our arguments.

\begin{lemma}(Dirac \cite{D}, Ore \cite{O})\label{le1}
Let $G$ be a graph of order $n\geq3$ and minimum degree $\delta(G).$ If  $\displaystyle \delta(G)\geq \frac{n+s}{2},$ then each of the following holds.\\
(i) $G$ is Hamiltonian-connected for $s=1$.\\
(ii) $G$ is Hamiltonian for $s=0$.\\
(iii) $G$ is traceable for $s=-1$.
\end{lemma}

\begin{lemma}(Chvat\'{a}l and Erd\"{o}s \cite{CE})\label{le2}
Let $G$ be a $k$-connected graph of order $n\geq3$. If $\displaystyle \alpha(G)\leq k-s,$ then each of the following holds.\\
(i) $G$ is Hamiltonian-connected for $s=1$.\\
(ii) $G$ is Hamiltonian for $s=0$.\\
(iii) $G$ is traceable for $s=-1$.
\end{lemma}

\begin{lemma}(Exercise 18.1.6, Bondy and Murty \cite{BM})\label{le4}
Let $G$ be a graph. Then $G$ is homogeneously traceable if and only if $K_{1}\nabla G$ is Hamiltonian-connected.
\end{lemma}


\begin{lemma}(Anderson and Morely \cite{AM})\label{le6}
Let $G$ be a graph of order $n\geq2$. Then $\mu_{1}(G)\leq n$ with equality if
and only if $\overline{G}$ is disconnected.
\end{lemma}

The main tool in our paper is the eigenvalue interlacing technique described below.
Given two non-increasing real sequences
$\theta_{1}\geq \theta_{2}\geq \cdots \geq \theta_{n}$ and
$\eta_{1}\geq \eta_{2}\geq \cdots \geq \eta_{m}$
with $n>m,$ the second sequence is said to {\bf interlace}
the first one if $\theta_{i}\geq \eta_{i}\geq\theta_{n-m+i}$
for $i=1, 2, \ldots, m.$
The interlacing is {\bf tight} if exists an integer $k\in[0, m]$
such that $\theta_{i}=\eta_{i}$ for $1\leq i\leq k$ and $\theta_{n-m+i}=\eta_{i}$ for
$k+1\leq i\leq m.$

Consider an $n\times n$ real symmetric matrix
\[
M=\left(\begin{array}{ccccccc}
M_{1,1}&M_{1,2}&\cdots &M_{1,m}\\
M_{2,1}&M_{2,2}&\cdots &M_{2,m}\\
\vdots& \vdots& \ddots& \vdots\\
M_{m,1}&M_{m,2}&\cdots &M_{m,m}\\
\end{array}\right),
\]
whose rows and columns are partitioned according to a partitioning
$X_{1}, X_{2},\ldots ,X_{m}$ of $\{1,2,\ldots, n\}$. The {\bf quotient matrix}
$R(M)$ of the matrix $M$ is the $m\times m$ matrix whose entries are the
average row sums of the blocks $M_{i,j}$ of $M$. The partition is {\bf equitable}
if each block $M_{i,j}$ of $M$ has constant row (and column) sum.

\begin{lemma}(Brouwer and Haemers \cite{BrHa12, Haem95})\label{le3}
Let $M$ be a real symmetric matrix. Then the eigenvalues of every
quotient matrix of $M$ interlace the ones of $M.$ Furthermore, if the
interlacing is tight, then the partition is equitable.
\end{lemma}

\section{Proofs}

\noindent{\bf Proof of Theorem \ref{main1}. } By observation, $K_{k,k}$ is not Hamiltonian-connected since there does not exist a Hamiltonian path between any pair of vertices belonging to the same part of the bipartition of $G.$ Therefore, it suffices to prove the sufficiency.

By contradiction, we assume that $G\ncong K_{k, k}$ and $G$ is not Hamiltonian-connected. First, we prove the following Claim.

\noindent {\bf Claim.} $n\geq 2k$.

In fact, if $n\leq 2k-1,$ then $\delta(G)\geq \kappa(G)\geq k\geq \frac{n+1}{2}$. By Lemma \ref{le1}(i), $G$ is Hamiltonian-connected, a contradiction.
Claim is completed.

By Lemma \ref{le2}(i), then $\alpha(G)\geq k,$ and thus there exists an independent set $X$ in $G$ such that $|X|=k.$ Let $X=\{u_{1}, u_{2}, \ldots, u_{k}\}$ and $Y=V(G)-X=\{v_{1}, v_{2}, \ldots, v_{n-k}\}$. Let $d_{1}(v_{i})=|N(v_{i})\cap X|$ and $d_{2}(v_{i})=|N(v_{i})\cap Y|$ for $1\leq i \leq n-k$. Clearly $\sum_{i=1}^{k}d(u_{i})=\sum_{i=1}^{n-k}d_{1}(v_{i})$.
Accordingly, the quotient matrix $R(A)$ of $A(G)$ on the partition $(X, Y)$ is as follows.
\[
R(A)=\left(\begin{array}{ccccccc}
0&\frac{\sum_{i=1}^{k}d(u_{i})}{k}\\
\frac{\sum_{i=1}^{n-k}d_{1}(v_{i})}{n-k}&\frac{\sum_{i=1}^{n-k}d_{2}(v_{i})}{n-k}\\
\end{array}\right)
=\left(\begin{array}{ccccccc}
0&\frac{\sum_{i=1}^{k}d(u_{i})}{k}\\
\frac{\sum_{i=1}^{k}d(u_{i})}{n-k}&\frac{\sum_{i=1}^{n-k}d_{2}(v_{i})}{n-k}\\
\end{array}\right).
\]
Let $\lambda_{1}(R(A))\geq \lambda_{2}(R(A))$ be the eigenvalues of $R(A)$.
By Lemma \ref{le3}, $$\lambda_{1}(G)\geq \lambda_{1}(R(A))\geq \lambda_{n-1}(G), ~~\lambda_{2}(G)\geq \lambda_{2}(R(A))\geq \lambda_{n}(G).$$
From Perron-Frobenius Theorem, we have $\lambda_{1}(G)\geq |\lambda_{n}(G)|$. Hence
\begin{eqnarray} \label{A-1}
\lambda_{1}^{2}(G) & \geq &-\lambda_{1}(G)\lambda_{n}(G) \geq -\lambda_{1}(R(A))\lambda_{2}(R(A))=-det(R(A))
\\ \nonumber
&=&\frac{\sum_{i=1}^{k}d(u_{i})}{k}\cdot \frac{\sum_{i=1}^{k}d(u_{i})}{n-k}=\frac{\sum_{i=1}^{k}d(u_{i})}{k}\cdot \frac{\sum_{i=1}^{k}d(u_{i})}{k}\cdot \frac{k}{n-k}
\\ \nonumber
&\geq& \frac{k\delta(G)}{k}\cdot \frac{k\delta(G)}{k}\cdot \frac{k}{n-k}=\frac{k\delta^{2}(G)}{n-k}
\geq \lambda_{1}^{2}(G).
\end{eqnarray}
It follows that all the inequalities in (\ref{A-1}) must be equalities. Hence we must have $\lambda_{1}(G)=-\lambda_{n}(G)$, $\lambda_{1}(G)=\lambda_{1}(R(A))$, $\lambda_{n}(G)=\lambda_{2}(R(A))$ and $d(u_{i})=\delta(G)$ for $1\leq i\leq k$. So $0=\lambda_{1}(G)+\lambda_{n}(G)=\lambda_{1}(R(A))+\lambda_{2}(R(A))=\frac{\sum_{i=1}^{n-k}d_{2}(v_{i})}{n-k}$, and thus $d_{2}(v_{i})=0$ for $1\leq i \leq n-k$. Therefore, $G$ is a bipartite graph with  the partition $(X, Y)$. Since
\begin{eqnarray*}
\delta(G)=\frac{\sum_{i=1}^{k}d(u_{i})}{k}=\frac{\sum_{i=1}^{n-k}d_{1}(v_{i})}{n-k}\cdot \frac{n-k}{k}=\frac{\sum_{i=1}^{n-k}d(v_{i})}{n-k}\cdot \frac{n-k}{k}\geq \delta(G) \cdot \frac{n-k}{k},
\end{eqnarray*}
we have $n\leq 2k.$ By Claim, then $n=2k$. Hence $n-k=k$ and $d(v_{i})=\delta(G)$. Note that $d(u_{i})=d(v_{i})=\delta(G)\geq \kappa(G)\geq k.$ Then $G\cong K_{k, k},$ contrary to our assumption. This completes the proof of the theorem.\hspace*{\fill}$\Box$

\noindent{\bf Proof of Corollary \ref{main2}.}~By assumption, $K_{1}\nabla G$ is $k+1$-connected of order $n+1$ and minimum degree $\delta(G)+1.$
Note that $K_{1}\nabla G\ncong K_{k, k}.$ By Theorem \ref{main1}, $K_{1}\nabla G$ is Hamiltonian-connected. By Lemma \ref{le4}, then $G$ is homogeneously traceable. \hspace*{\fill}$\Box$

In the proof of Theorem \ref{main3}, we say that a graph possesses Hamiltonian properties, which means that the graph is Hamiltonian-connected, Hamiltonian or traceable.

\noindent{\bf Proof of Theorem \ref{main3}.}~ We divide our proof into two cases.

\noindent {\bf Case 1.} $\delta > n-k+s-1$.

As $\delta \geq k$ and $\delta > n-k+s-1$, we have $\delta > \frac{n+s-1}{2}.$

\noindent {\bf Case 1.1.} $n=2t+1$ is odd, where $t\geq 1$.

If $s=1$, then $\delta > \frac{2t+1}{2}$, and we have $\delta \geq t+1=\frac{n+s}{2}.$ By Lemma \ref{le1}(i), $G$ is Hamiltonian-connected.
If $s=0$, then $\delta > t$, and we have $\delta \geq t+1>\frac{n+s}{2}.$ By Lemma \ref{le1}(ii), $G$ is Hamiltonian.
If $s=-1$, then $\delta > \frac{2t-1}{2}$, and we have $\delta \geq t=\frac{n+s}{2}.$ By Lemma \ref{le1}(iii), $G$ is traceable.

\noindent
{\bf Case 1.2.} $n=2t+2$ is even, where $t\geq 1$.

If $s=1$, then $\delta > t+1$, and we have $\delta \geq t+2>\frac{n+s}{2}.$ By Lemma \ref{le1}(i), $G$ is Hamiltonian-connected.
If $s=0$, then $\delta > t+\frac{1}{2}$, and we have $\delta \geq t+1=\frac{n+s}{2}.$ By Lemma \ref{le1}(ii), $G$ is Hamiltonian.
If $s=-1$, then $\delta > t$, and we have $\delta \geq t+1>\frac{n+s}{2}.$ By Lemma \ref{le1}(iii), $G$ is traceable.

\noindent {\bf Case 2.} $\delta\leq n-k+s-1$.

Now, $\frac{n\delta}{n-k+s-1}\leq n.$ Suppose, to the contrary, that $G$ does not possess Hamiltonian properties.
By Lemma \ref{le2}, $\alpha(G)\geq k-s+1,$ and then there exists an independent set $X$ in $G$ such that $|X|=k-s+1.$ Let $X=\{u_{1}, u_{2}, \ldots, u_{k-s+1}\}$ and $Y=V(G)-X=\{v_{1}, v_{2}, \ldots, v_{n-k+s-1}\}$. Let $d_{1}(v_{i})=|N(v_{i})\cap X|$ and $d_{2}(v_{i})=|N(v_{i})\cap Y|$ for $1\leq i \leq n-k+s-1$. Clearly $\sum_{i=1}^{k-s+1}d(u_{i})=\sum_{i=1}^{n-k+s-1}d_{1}(v_{i})$.
Accordingly, the quotient matrix $R(L)$ of $L(G)$ on the partition $(X, Y)$ becomes:
\[
R(L)=\left(\begin{array}{ccccccc}
\frac{\sum_{i=1}^{k-s+1}d(u_{i})}{k-s+1}&-\frac{\sum_{i=1}^{k-s+1}d(u_{i})}{k-s+1}\\
-\frac{\sum_{i=1}^{n-k+s-1}d_{1}(v_{i})}{n-k+s-1}&\frac{\sum_{i=1}^{n-k+s-1}d_{1}(v_{i})}{n-k+s-1}\\
\end{array}\right)
=\left(\begin{array}{ccccccc}
\frac{\sum_{i=1}^{k-s+1}d(u_{i})}{k-s+1}&-\frac{\sum_{i=1}^{k-s+1}d(u_{i})}{k-s+1}\\
-\frac{\sum_{i=1}^{k-s+1}d(u_{i})}{n-k+s-1}&\frac{\sum_{i=1}^{k-s+1}d(u_{i})}{n-k+s-1}\\
\end{array}\right).
\]
Let $\mu_{1}(R(L))\geq \mu_{2}(R(L))=0$ be the eigenvalues of $R(L)$. By algebraic manipulation, we have $$\mu_{1}(R(L))=(\frac{1}{k-s+1}+\frac{1}{n-k+s-1})\sum_{i=1}^{k-s+1}d(u_{i}).$$ By Lemma \ref{le3}, then
\begin{eqnarray*}
\mu_{1}(G)& \geq & \mu_{1}(R(L))=(\frac{1}{k-s+1}+\frac{1}{n-k+s-1})\sum_{i=1}^{k-s+1}d(u_{i})
\\ \nonumber
& \geq & (\frac{1}{k-s+1}+\frac{1}{n-k+s-1})(k-s+1)\delta
= \frac{n\delta}{n-k+s-1}.
\end{eqnarray*}
This contradicts the assumption of this theorem.\hspace*{\fill}$\Box$

\vspace{5mm}
\noindent
{\bf Acknowledgement.} The research of Ruifang Liu is partially
supported by NSFC (No.~11571323) and Foundation for University Key Teacher of Henan Province (No.~2016GGJS-007).

\end{document}